\documentclass{article}

\usepackage{arxiv}

\usepackage[utf8]{inputenc} 
\usepackage[T1]{fontenc}    
\usepackage{hyperref}
\hypersetup{colorlinks=true,linkcolor=red,citecolor=blue}
\usepackage{url}            
\usepackage{booktabs}       
\usepackage{amsfonts}       
\usepackage{nicefrac}       
\usepackage{microtype}      
\usepackage{graphicx}
\usepackage[square,numbers]{natbib}   
\usepackage{subfigure}
\usepackage{xcolor}
\usepackage{array}
\usepackage{xpatch}
\usepackage{amsmath,amssymb,amscd,bm}
\usepackage{mathrsfs}
\usepackage{paralist}
\usepackage{tikz}
\usetikzlibrary{3d}
\usetikzlibrary{calc}

\title{
Physics-Aware Neural Networks for Boundary Layer Linear Problems
}

\author{
Ant\^onio Tadeu Azevedo Gomes \\
Laborat\'orio Nacional de Computa\c c\~ao Cient\'ifica (LNCC) \\
\texttt{atagomes@lncc.br}
\And
Larissa Miguez da Silva \\
Programa de Doutorado em Modelagem Computacional -- Laborat\'orio Nacional de Computa\c c\~ao Cient\'ifica (LNCC) \\
\texttt{lamiguez@posgrad.lncc.br} \\
\And
Fr\'ed\'eric Valentin \\
Laborat\'orio Nacional de Computa\c c\~ao Cient\'ifica (LNCC) \\
\texttt{valentin@lncc.br}
}

\date{}

\renewcommand{\headeright}{}
\renewcommand{\undertitle}{}

\hypersetup{
}
\begin{document}
\maketitle              
\begin{abstract}
\textit{Physics-Informed Neural Networks} (PINNs) are machine learning tools that approximate the solution of general partial differential equations (PDEs) by adding them in some form as terms of the loss/cost function of a Neural Network.
Most pieces of work in the area of PINNs tackle \textit{non-linear} PDEs.
Nevertheless, many interesting problems involving linear PDEs may benefit from PINNs; these include parametric studies, multi-query problems, and parabolic (transient) PDEs.
The purpose of this paper is to explore PINNs for linear PDEs whose solutions may present one or more \textit{boundary layers}.
More specifically, we analyze the steady-state reaction-advection-diffusion equation in regimes in which the diffusive coefficient is small in comparison with the  reactive or advective coefficients.
We show that adding information about these coefficients as predictor variables in a PINN results in better prediction models than a PINN that only uses spatial information as predictor variables.
This finding may be instrumental in multiscale problems where the coefficients of the PDEs present high variability in small spatiotemporal regions of the domain, and therefore PINNs may be employed together with domain decomposition techniques to efficiently approximate the PDEs locally at each partition of the spatiotemporal domain, without resorting to different learned PINN models at each of these partitions.

\end{abstract}

\section{Introduction}
\textit{Physics-Aware} Neural Networks (NNs) are machine learning tools that approximate the solution of general partial differential equations (PDEs) by adding the physical laws these equations represent to some component of the neural network. 
The PINNs~\cite{RaiPerKar2019} are likely to be the most well-known of these NNs; the defining characteristic of a PINN is the inclusion of the strong form of a PDE (including its boundary and initial conditions) as terms of the loss/cost function. 

Most pieces of work in the area of physics-aware NNs tackle \textit{non-linear} PDEs~\cite{EYu2018,KhaZhaKar2019,LiaMin2021,JagKar2020,JagKhaKar2020}.
Nevertheless, many interesting problems involving linear PDEs may benefit from physics-aware NNs; these include parametric studies, multi-query problems, and parabolic (transient) PDEs.

We are mostly interested in the solution of linear PDEs whose coefficients present high variability in small spatiotemporal regions of the physical domain.
In this case, we say that the solution has a multiscale behavior. 
Standard numerical methods often present difficulties in approximating the solution to such PDEs with combined quality and computational affordability. 
Multiscale numerical methods (e.g.~\cite{harder2015multiscale}) have emerged as an attractive option for dealing with such difficulties by rewriting the original formulation of the PDE in terms of:
\begin{inparaenum}[(i)]
\item local problems living each one in a partition of the physical domain; and
\item a global problem that ``glues together'' the solution of the local problems.
\end{inparaenum}
The price to pay is a potentially large number of local problems. 
Although said local problems are independent from one another, thus benefiting from massive parallel computations, they may still be computationally demanding.

The purpose of this paper is to investigate the potential of physics-aware NNs in general, and PINNs specifically, for efficiently solving local problems in multiscale numerical methods.
We explore the particular case of linear PDEs whose solutions may present one or more \textit{boundary layers}.
More specifically, we analyze the steady-state reaction-advection-diffusion equation in regimes in which the diffusive coefficient is small in comparison with the reactive or advective coefficients.
We verify that adding information about these coefficients as predictor variables in a PINN results in better prediction models than a PINN that only uses spatial information as predictor variables.
We believe this finding may be instrumental in multiscale problems, because it opens the path for PINNs to be employed together with domain decomposition techniques to efficiently approximate the PDEs locally at each partition of the spatiotemporal domain, without resorting to different learned PINN models at each of these partitions.

The remainder of this paper is structured as follows.
In Section~\ref{sec:relwork}, we quickly review the related literature. 
In Section~\ref{sec:methodology}, we present the problem and the methodology for the proposed model. 
In Section~\ref{sec:results}, we examine two different cases of the target equation and the effectiveness of the proposed model. 
Finally, in Section~\ref{sec:conclusion}, we report the conclusions of this work along with a discussion of future directions.

\section{Related work}
\label{sec:relwork}

In recent years, the use of algorithms that ``learn'' from data has caused great impact and change in several areas of science. 
Algorithms using NNs have been used in many problems governed by PDEs and presented satisfactory results~\cite{berg2018unified, lagaris1998artificial, raissi2018hidden, sirignano2018dgm}. 
In particular, in \cite{raissi2019physics} the methodology of PINNs was first proposed, combining the properties of universal approximation of NNs and the knowledge of physical laws described by PDEs. 
Since then, many pieces of work have been published on this topic~\cite{mao2020physics, cai2021physics, misyris2020physics}. 
However, problems with complex geometry domains have led to other methodologies based on domain decomposition methods and PINNs,  including Extended PINNs (XPINNs)~\cite{jagtap2021extended}, Conservative PINNs (cPINNs)~\cite{jagtap2020conservative} and Variational PINNs (VPINNs)~\cite{kharazmi2019variational}.

Although the aforementioned pieces of work have shown excellent results, there are still many theoretical gaps that need to be filled. 
The techniques are new and do not have a trivial application in the solution of physical problems. 
For instance, a particularly complex step in formulating deep learning problems and PINNs is the definition of the loss functional to be minimized. 
Additionally, there are many hyperparameters to be configured and, although the automatic selection of hyperparameters is possible, there is usually a large computational cost. 
Interestingly, this last problem also exists somehow in a handful of multiscale methods with regard to their configuration parameters (e.g.~\cite{wscad2020,10.1007/978-3-030-68035-0_7}). 

There is a growing number of papers relating multiscale methods and data-driven approaches~\cite{capuano2019smart,chan2018machine,karpatne2017physics,yeung2020deep}. 
To the best of our knowledge, none of these pieces of work tackle the problem the way we do, which is by training a single machine learning model that may be parameterized for approximating the solution of a PDE with highly variable coefficients in the spatial domain.

\section{Methodology}
\label{sec:methodology}

\subsection{Boundary layer problem}
\label{subsec:physicalProblem}

The boundary layer problem can appear in many applications, including fluid dynamics, meteorology, atomic reactors, among others. This phenomenon occurs when the gradient is high in the region close to the boundary and can bring instability to the discrete solution of the problem. Next, we present an example for this case.

Consider the case in which the reaction-advection-diffusion problem has an exact solution which contains boundary layers. 
This happens when the reactive or advective coefficient dominates the diffusive one. 
We consider the following reaction-advection-diffusion problem: \textit{Find} $u \in H^1(\Omega)$ \textit{such that}:
\begin{equation}
    \begin{cases}
    \begin{aligned}
\nabla \cdot (- \mathcal{K} \nabla u + \alpha u) + \sigma u &= 1 \text{  in }  \Omega,\\
u &= 0 \text{  on } \partial \Omega_D, \\
\nabla u \cdot \mathbf{n} &= 0 \text{  on } \partial \Omega_N,
\end{aligned}
\end{cases}\label{EQ:ARD}
\end{equation}
where $\Omega$ is a unit square domain, $\partial \Omega_D$ corresponds to the boundaries $x = (0, y)$ and $x = (1, y)$, with $y \in (0, 1)$, where homogeneous Dirichlet conditions are to be enforced, and $\partial \Omega_N = \partial \Omega \backslash \partial \Omega_D$ corresponds to the boundaries where homogeneous Neumann conditions are to be enforced.
The coefficients are such that $\alpha := (a, 0)^T$, $\mathcal{K} =  k \mathcal{I} $, where $\mathcal{I}$ is the identity matrix and $a$, $k$, $\sigma \in \mathbb{R}$. If $\sigma>0$, the analytical solution to this equation is:

\begin{equation*}
    u(x,y) = \dfrac{\sinh( \dfrac{\sqrt{4 k \sigma}}{2 k}(x-1)) - \sinh(\dfrac{\sqrt{4 k \sigma}}{2 k} x)}{\sinh(\dfrac{\sqrt{a^2 + 4 k \sigma}}{2 k})} +1\,,
\end{equation*}

Otherwise, if $\sigma=0$ and $a > 0$, then the exact solution becomes 
\begin{equation*}
    u(x,y) = \dfrac{1}{a} \left( x - \dfrac{\sinh( \dfrac{\sqrt{a}}{2 k}x)}{\sinh(\dfrac{\sqrt{a}}{2 k})} e^{\frac{\sqrt{a}}{2 k}(x-1) }\right)
\end{equation*}

We consider two experimental settings.
First, in Subsection~\ref{sec:results-reactive} we will assume $\sigma=1$ and $a=0$. 
So, we will explore the case in which the reactive coefficient dominates the diffusive one. 
Next, in Subsection~\ref{sec:results-advective} we will consider the case in which the advective coefficient dominates by assuming that $\sigma=0$ and $a=1$.

\subsection{Architecture design}
\label{subsec:architecture}
\begin{figure}[h!]
    \centering
    \scalebox{0.8}{\input{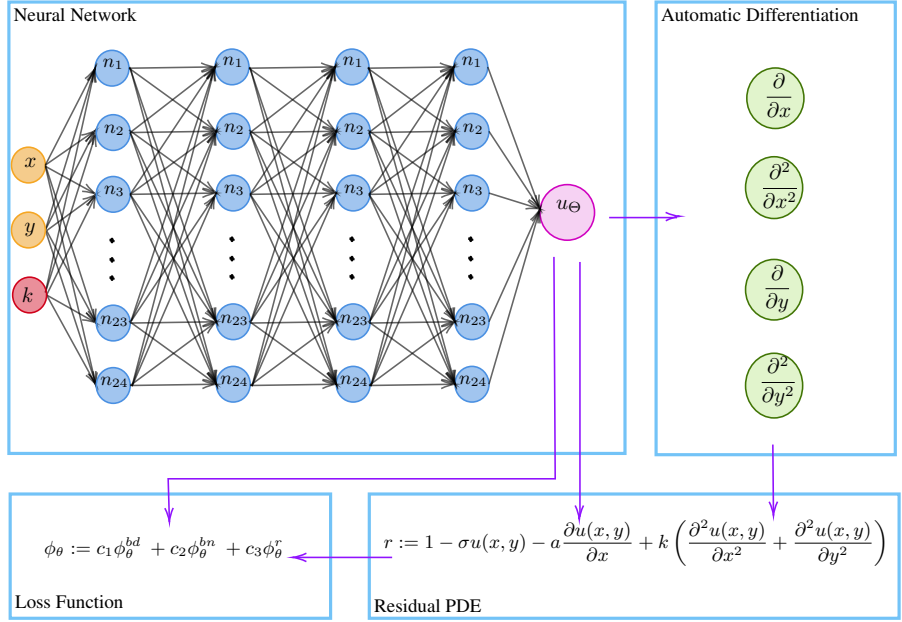}}
    \caption{PINN architecture for Reaction-Advection-Diffusion problem.}
    \label{fig:pinns_architecture}
\end{figure}

To solve the problem presented in Subsection~\ref{subsec:physicalProblem}, we use the PINN depicted in Figure~\ref{fig:pinns_architecture}.
A key feature of this PINN is the use of the diffusion coefficient $k$ as a predictor variable together with the spatial data $(x,y)$.
This way, we aim at getting a model capable of making predictions for different diffusion coefficients.

We assume a feed-forward NN with the following structure: 4 fully connected layers each containing 24 neurons\footnote{These hyperparameters as well as all other configurations not explicitly explained in the remainder of the text have been determined empirically.} and each followed by a hyperbolic tangent activation function. 
Furthermore, we use one output layer of size 1 and a linear activation function. 
For the sake of comparison, we establish two scenarios for the input layer, as will be further explained in the following section:
\begin{inparaenum}[(i)]
\item \textit{Scenario 1}, with only 2 neurons, input $k$ being taken out;
\item \textit{Scenario 2} with 3 neurons, exactly as shown in Figure~\ref{fig:pinns_architecture}.
\end{inparaenum}

Also, we considered the loss function 

\begin{equation}
    \phi_\theta(X) :=  c_1 \phi_\theta^{bd}(X^{bd}) + c_2 \phi_\theta^{bn}(X^{bn}) + c_3 \phi_\theta^r(X^r),
\end{equation}

\noindent as a function of the training data, as also explained in the following section.

\section{Experimental results}
\label{sec:results}

In this section, we present some numerical results that show the performance of PINNs to solve the boundary layer problem. 
In our simulations, we consider the following scenarios:
\begin{itemize}
    \item[\textit{Scenario 1}:] We fix a diffusion coefficient $k$, train the network for some collocation and boundary points and predict others. Therefore, the input parameters of the network are $(x,y)$ and the output is the solution $u$.
    \item[\textit{Scenario 2}:] We vary the diffusion coefficient $k$ and train the network for some $k$ to predict. Therefore, the input parameters of the network are $(x, y, k)$ and the output is the solution $u$.
\end{itemize}

The collocation points are given by $X_r$, the boundary data is in $X_{bd}$ and $X_{bn}$, where $X_{bd}$ represents the data on the Dirichlet boundary $\Omega_D$ and $X_{bn}$ the data on the Neumann boundary $\Omega_N$. 
Respectively, on those boundary points, we have the solutions $u_{bd}$ and $u_{bn}$. 
Additionally, the coefficients $c_1$, $c_2$ and $c_3$ representing the weights of each loss term are hyperparameters of the model, and their values have been chosen empirically based on the knowledge of the authors about the behavior of Equation~\ref{EQ:ARD} for different values of its coefficients.

We assume that the collocation points $X_r$ as well as the points for the boundary data $X_{bd}$ and $X_{bn}$ are randomly sampled from a uniform distribution. 

\subsection{First setting: reaction-diffusion problem}
\label{sec:results-reactive}

\begin{figure}[h!]
    \centering
    \includegraphics[width=0.60\textwidth]{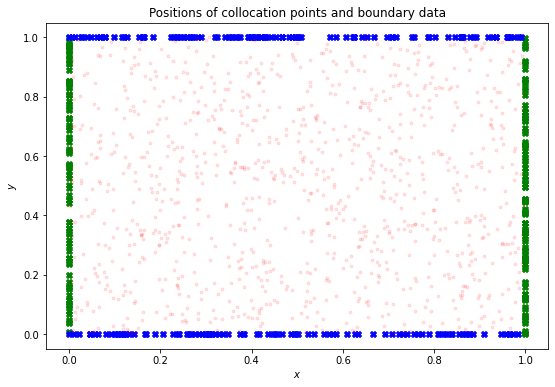}
    \caption{The collocation points (red circles) and the positions where the Dirichlet boundary condition (blue cross marks) and the Neumann boundary condition (green cross marks) will be weakly imposed.}
    \label{fig:points_of_collocation_and_boundary_data}
\end{figure}

For this first problem, we began with a training data of size $N_{bd}=N_{bn}=200$ and $N_r=1000$, where $N_{bd}$ is when we apply the Dirichlet boundary condition and $N_{bn}$ when we apply the Neumann boundary condition. 
We illustrate this setting in Figure~\ref{fig:points_of_collocation_and_boundary_data}.

\begin{figure}[h!]
    \centering
    \subfigure[Predicted solution ($k = 0.0001$) ]{
        \includegraphics[width=.40\textwidth]{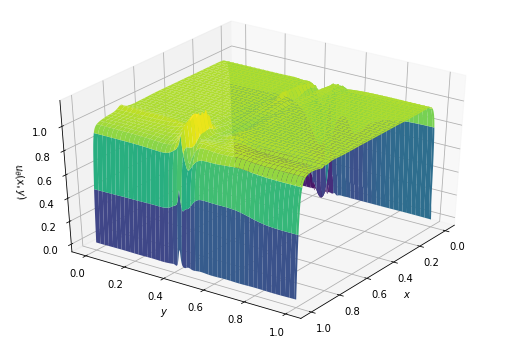}
    } 
    \quad 
    \subfigure[Exact solution ($k = 0.0001$) ]{
        \includegraphics[width=.40\textwidth]{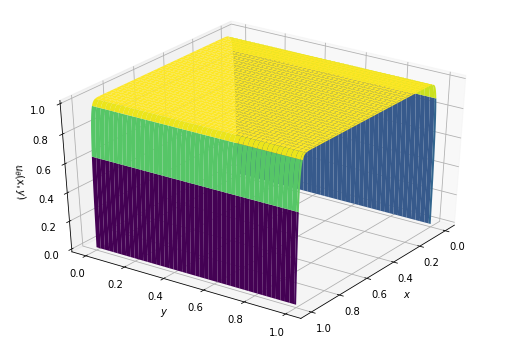}
        } 
        \\
    \subfigure[Predicted solution ($k = 0.00001$) ]{
    \includegraphics[width=.40\textwidth]{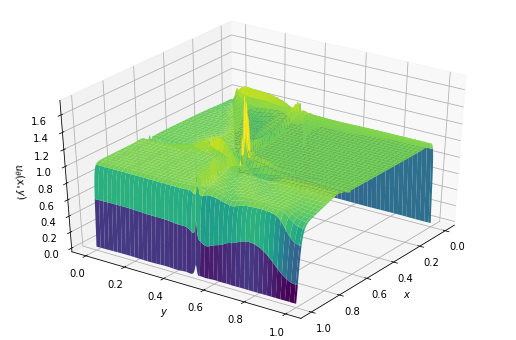}
    }
    \quad 
    \subfigure[Exact solution ($k = 0.00001$)]{
    \includegraphics[width=.40\textwidth]{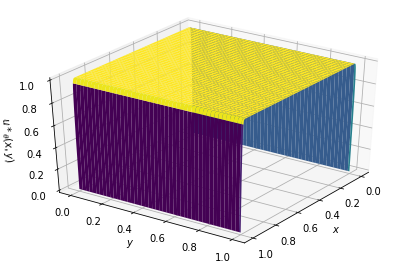}
    }
\caption{Scenario 1, reaction-diffusion setting: predicted solution vs exact solution with respect to parameter $k$.}
 \label{fig:case1_results}
\end{figure}

\begin{figure}[h!]
    \centering
    \includegraphics[width=.5\textwidth]{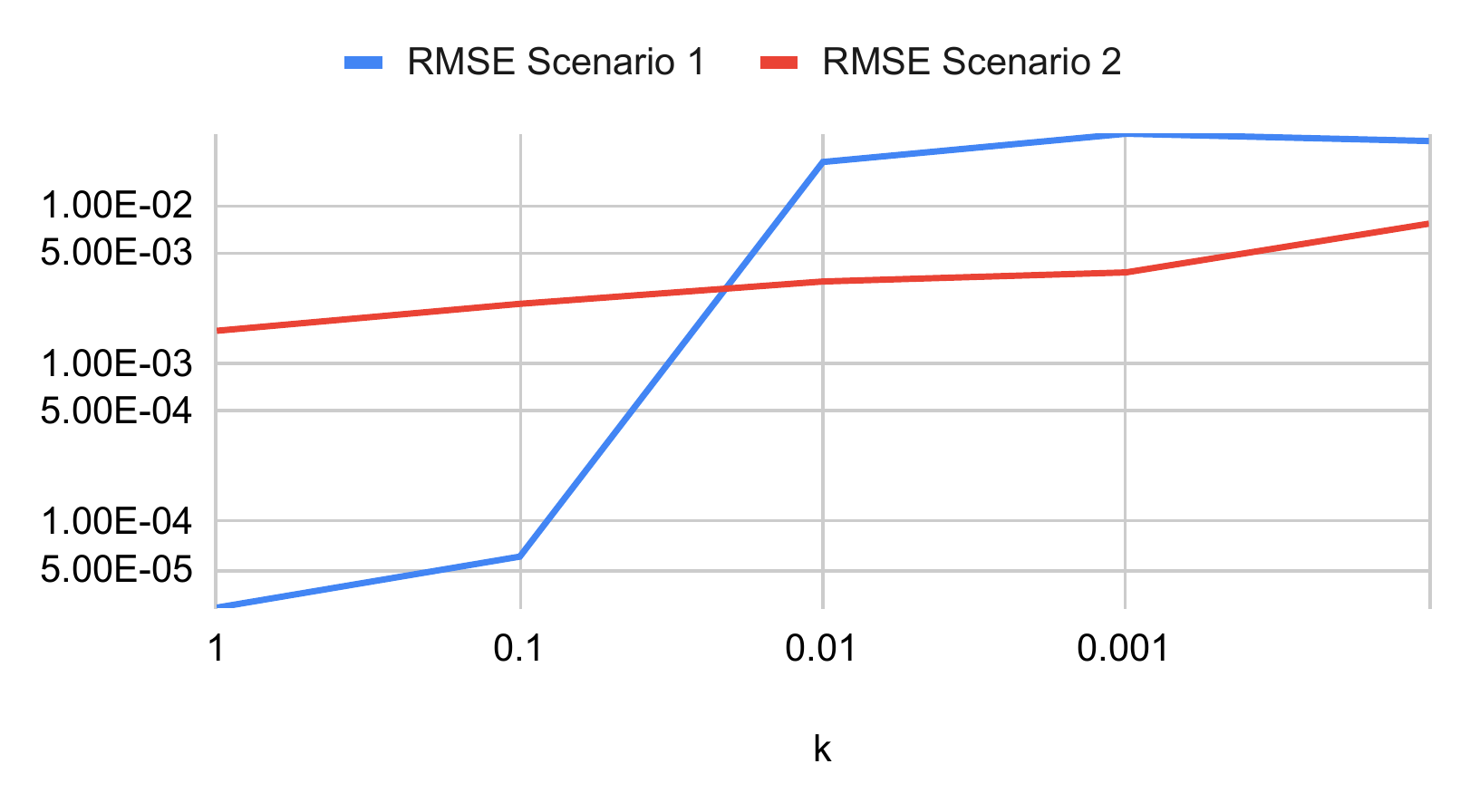}
    \caption{Quality of prediction for the reaction dominant case (log-log scale).}
    \label{fig:reaction_case}
\end{figure}

Figure~\ref{fig:case1_results} depicts some experimental results for Scenario 1.  
We observe that as we shrink the diffusive coefficient $k$, the solution gets worse, with some undesired features when $k = 0.0001$ and $k = 0.00001$. 
(Nonetheless, we see in Figure~\ref{fig:reaction_case}
that the PINN is able to approximate the solution with a small error when $k = 1.0$, $k = 0.1$, and $k = 0.01$.)
Besides the difficulty of approximating the solution in the case where we have a very small $k$, another disadvantage of this approach is that the model needs to be retrained for each new $k$. 

\begin{figure}[h!]
    \centering
    \subfigure[$N_{bd}=N_{bn}=100$, $N_r=480$ ]{
        \includegraphics[width=.47\textwidth]{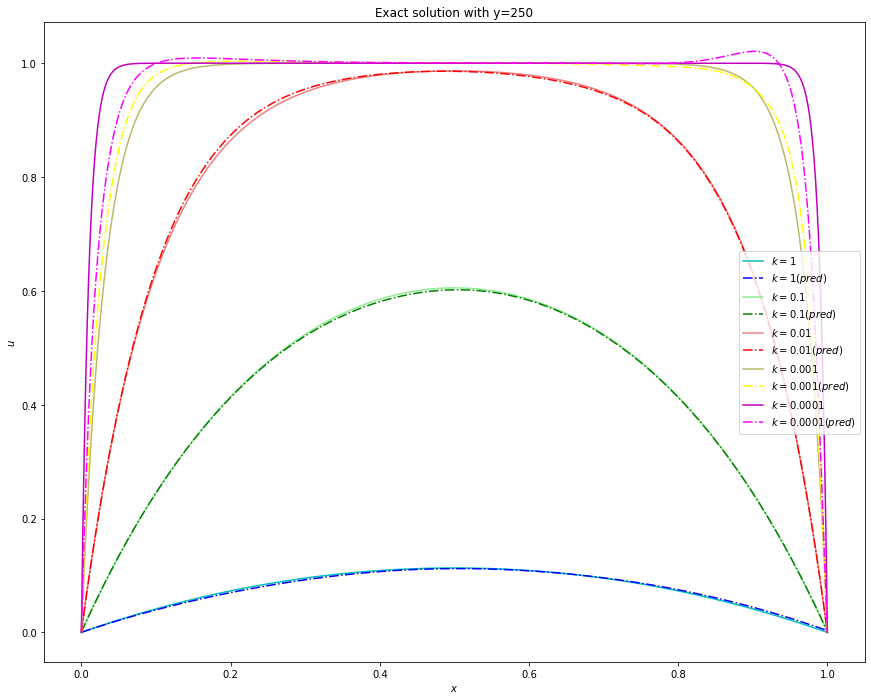}
    } 
    \subfigure[$N_{bd}=N_{bn}=150$, $N_r=600$ ]{
        \includegraphics[width=.47\textwidth]{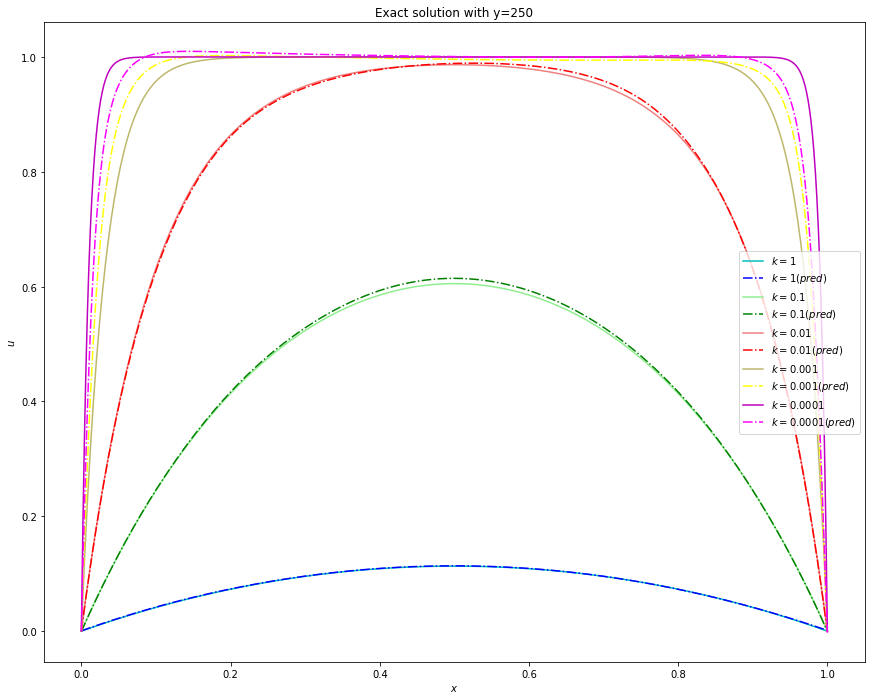}
        } 
        \\
         \subfigure[$N_{bd}=N_{bn}=200$, $N_r=800$ ]{
        \includegraphics[width=.47\textwidth]{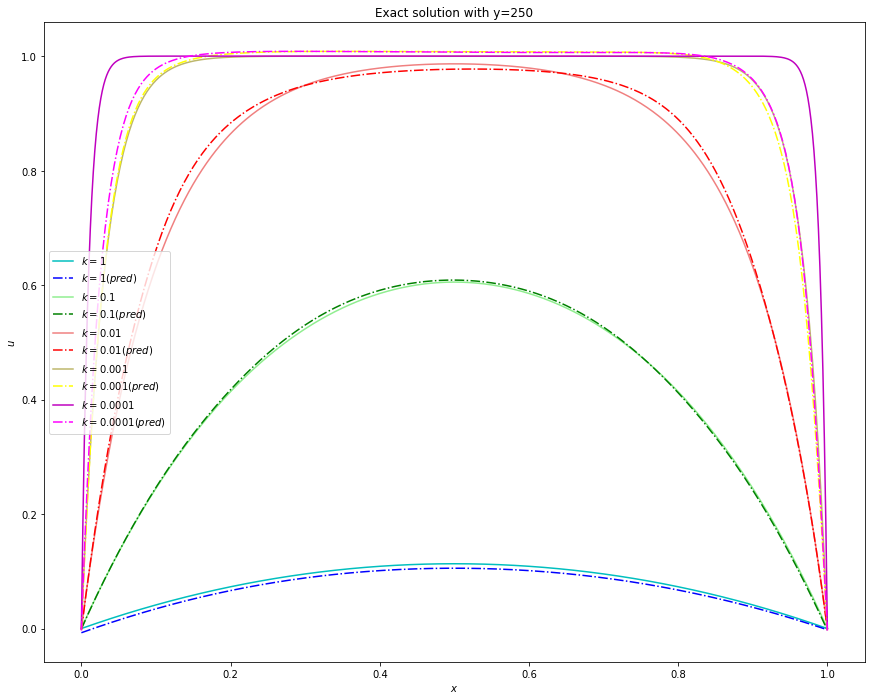}
        }
    \subfigure[$N_{bd}=N_{bn}=200$, $N_r=1000$]{
        \includegraphics[width=.47\textwidth]{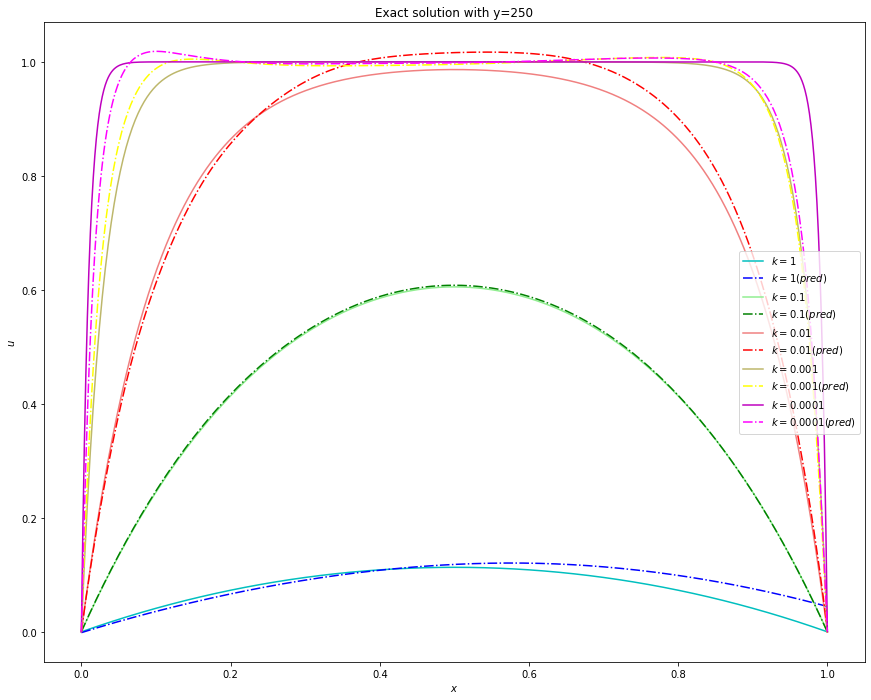}
        }
\caption{Scenario 2, reaction-diffusion setting: sensitivity analysis with respect to the size of the training data.}
 \label{fig:case2_collocationTrainingPoins}
\end{figure}

For Scenario 2, we add 20 different, randomly sampled $k \in (0.0001, 1.0)$ to the input data.
First, we investigate the sensitivity of the PINN with respect to the amount and dispersion of collocation training points in Scenario 2.
In Figure~\ref{fig:case2_collocationTrainingPoins}, we show the exact solution and the solution field generated by a PINN with decreasing values of $k$, for different amounts of collocation and boundary points.
We plot a cut for a fixed $y$; the exact solutions are represented by the solid lines and the predicted PINN solutions are represented by the dotted lines.

We can observe that the proposed PINN architecture for Scenario 2 interpolates quite well for values of $k$ greater than or equal to $0.001$, but for $k = 0.0001$ we have significant errors. 
Once more we see the impact of the boundary layer problem on the predictions. 

\begin{figure}[h!]
    \centering
    \subfigure[$c_3 = 0.2$ ]{
        \includegraphics[width=.35\textwidth]{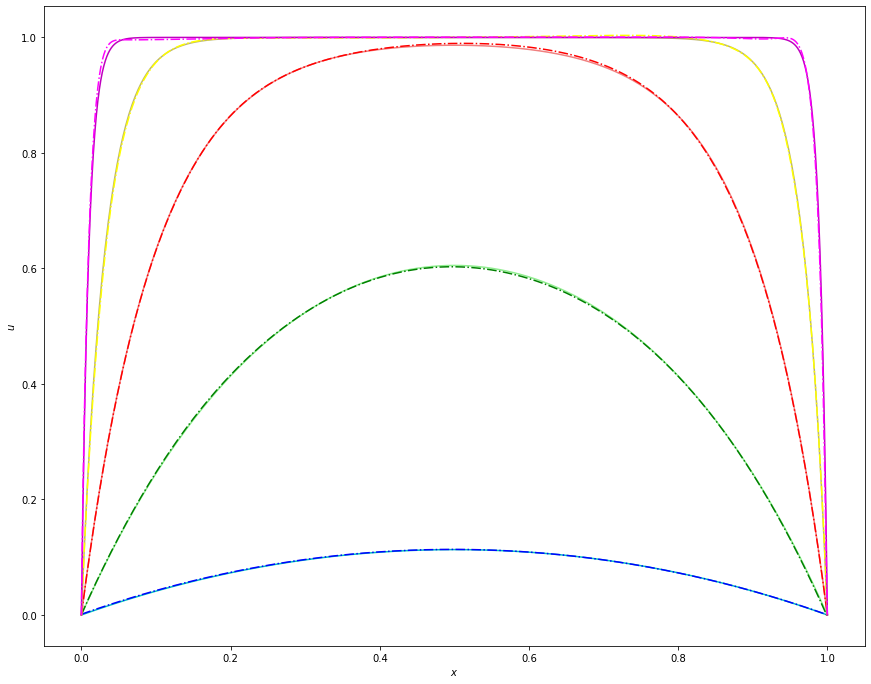}
    } 
    \subfigure[$c_3 = 0.4$ ]{
        \includegraphics[width=.35\textwidth]{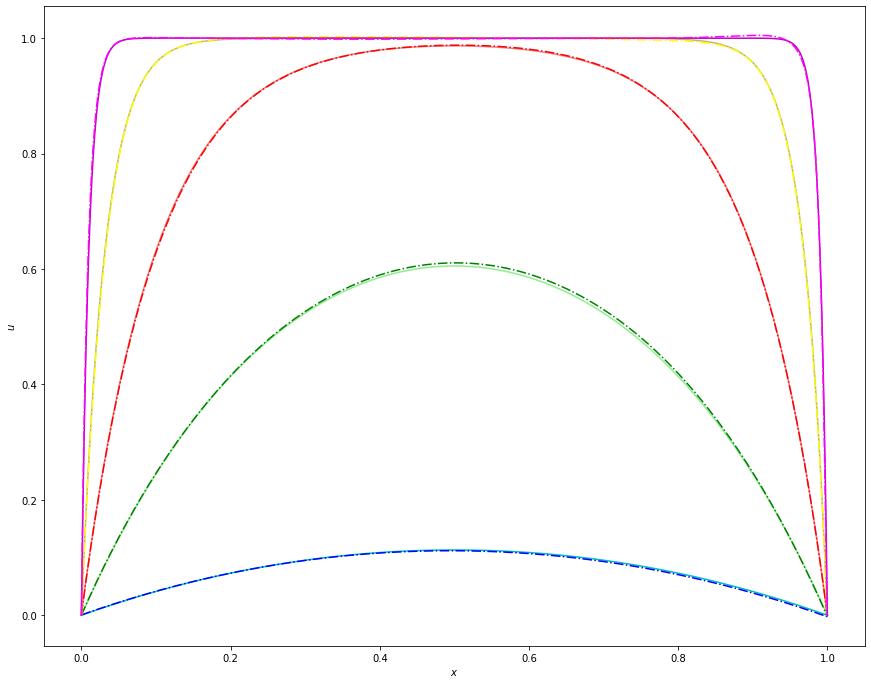}
        } 
    \includegraphics[width=.13\textwidth]{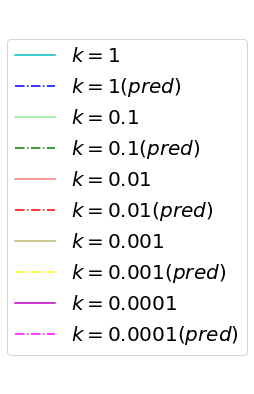}
\caption{Scenario 2, reaction-diffusion setting: Sensitivity analysis with respect to loss weights.}
    \label{fig:case2_lossWeights}
\end{figure}

In the above experiments we trained the PINNs with $c_1 = 2 $, $c_2 = 1 $ and $c_3 = 0.01$ as the weights of the loss terms. 
In Figure~\ref{fig:case2_lossWeights} we show the sensitivity of the model for Scenario 2 with respect to the loss weights.

For the cases presented in Figure \ref{fig:case2_lossWeights} we used $N_{bd}=N_{bn}=100$, $N_r=240$. 
The results are even more impressive because the PINN algorithm in Scenario 2 is able to reconstruct the solution field with high precision from a small number of points used for the training, even for the case where we have a very small $k$. 
Therefore, the results clearly show the impact of the loss weights on the training.

Finally, we compare the errors originating from Scenario 1 and Scenario 2. 
We use the Relative Mean Square Error (RMSE) to compare the results of the different scenarios, as presented in Figure~\ref{fig:reaction_case}.
We observe a significant error increase with a decaying $k$ for Scenario 1, whereas for Scenario 2 this increase is much slower, specially for $k \in (0.001, 1.0)$. 
However, the much higher errors for Scenario 2 with larger values of $k$ are still largely unexplained and motivates a series of investigations as part of our future work.

\begin{figure}[h!]
    \centering
    \subfigure[Extrapolating to $k = 10^{-5}$ ]{
        \includegraphics[width=.45\textwidth]{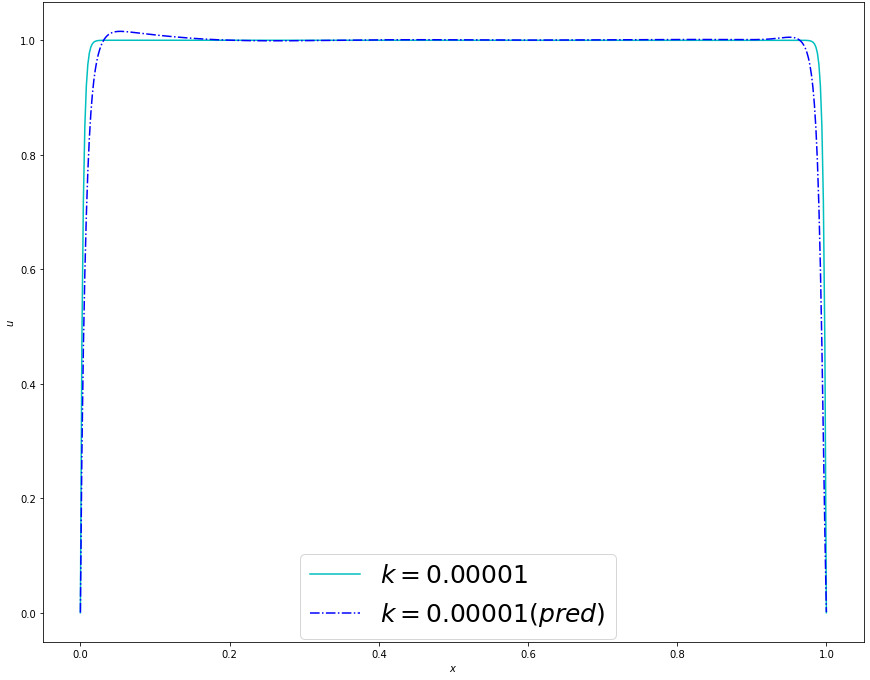}
    } 
    \quad 
    \subfigure[Extrapolating to $k = 1.2$ ]{
        \includegraphics[width=.45\textwidth]{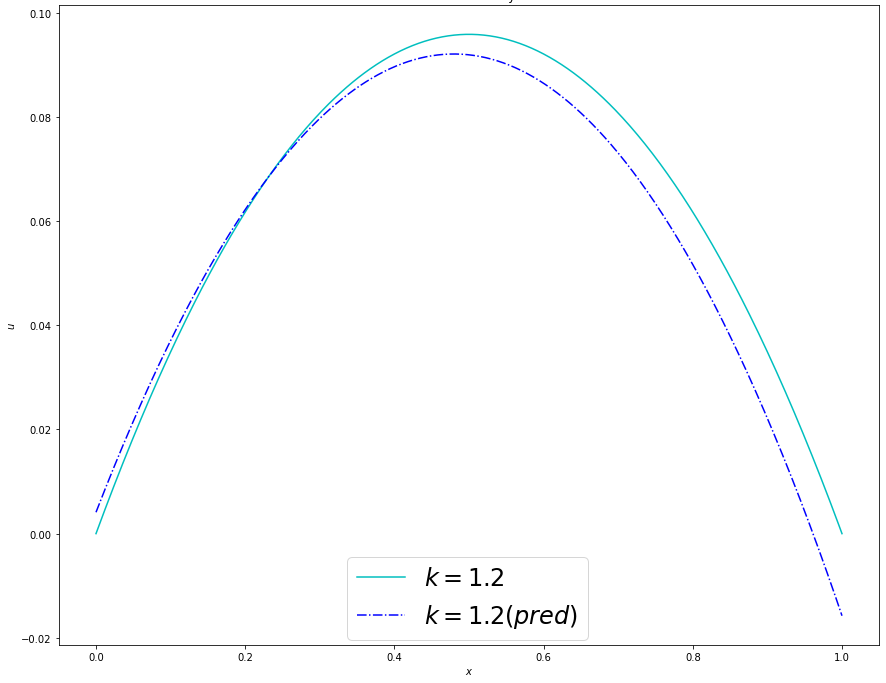}
        } 
\caption{Scenario 2, reaction-diffusion setting: extrapolating for different values of $k$.}
    \label{fig:case2_extrap}
\end{figure}

Now, we take advantage of the fact that we can predict for different values of $k$ in Scenario 2 and try to perform an extrapolation. 
The results are presented in Figure~\ref{fig:case2_extrap}, in which we plot two different cases. In Figure~\ref{fig:case2_extrap}(a), we extrapolate for an even smaller $k$, where the boundary layer problem is more evident. 
Even so, we still get a good approximation; we believe this is due to the fact that the diffusion coefficient is close to the range of $k$ used for the training. 
The same does not occur in Figure~\ref{fig:case2_extrap}(b), when we try to extrapolate to a larger value of $k$. 
(Mind, however, that we have plotted Figure~\ref{fig:case2_extrap}(b) in a different scale, to emphasize the prediction error.)
The curve for the exact solution is far from the curve that represents the predicted solution although this case is a completely boundary layer free problem. 
This result is likely an indication that the proposed method still does not work well for extrapolations far from the set used for the training.
Also note that the error is particularly high near the Dirichlet boundary $(1, y)$, with $y \in (0, 1)$, which shows that imposing the boundary condition weakly by means of a loss term may be tricky.

\subsection{Second setting: advection-diffusion problem}
\label{sec:results-advective}

Here, similarly to Subsection~\ref{sec:results-reactive}, we will consider the same two scenarios, as well as the same PINN architecture, and the same number of collocation points already described. 
What we will change are the loss weights, now set to $c_1 = 1$, $c_2 = 1.2$ and $c_3 = 1$. 
    
\begin{figure}[h!]
    \centering
    \subfigure[Predicted solution ($k = 0.01$) ]{
        \includegraphics[width=.42\textwidth]{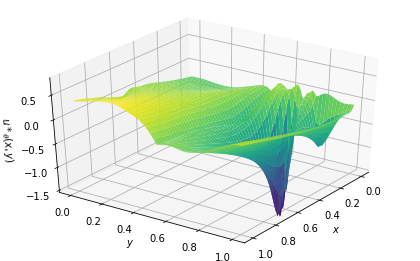}
    } 
    \quad 
    \subfigure[Exact solution ($k = 0.01$) ]{
        \includegraphics[width=.42\textwidth]{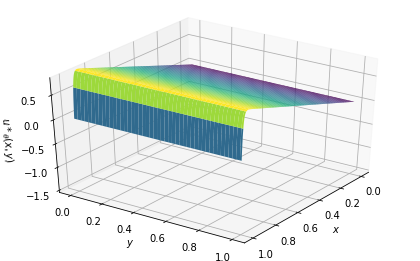}
        } 
        \\
    \subfigure[Predicted solution ($k = 0.001$) ]{
    \includegraphics[width=.42\textwidth]{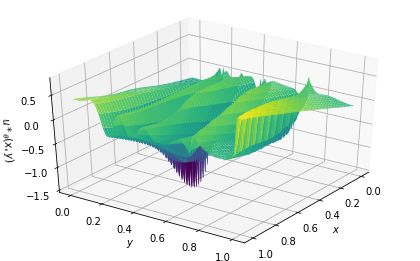}
    }
    \quad 
    \subfigure[Exact solution ($k = 0.001$)]{
    \includegraphics[width=.42\textwidth]{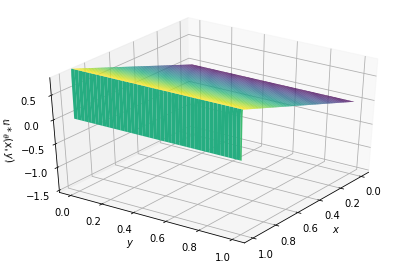}
    }
\caption{Scenario 1, advection-diffusion setting: predicted solution vs exact solution with respect to parameter $k$.}
 \label{fig:adv_results_scenario1}
\end{figure}

Figure~\ref{fig:adv_results_scenario1} depicts some experimental results for Scenario 1.  
As in the reaction-diffusion setting, we observe that as we shrink the diffusive coefficient $k$, the solution gets worse.
The advection-diffusion problem is nevertheless much tougher to approximate well than the reaction-diffusion problem.
The figure clearly shows this, with completely wrong solutions when $k = 0.01$ and $k = 0.001$. 

\begin{figure}[h!]
    \centering
    \includegraphics[width=0.45\textwidth]{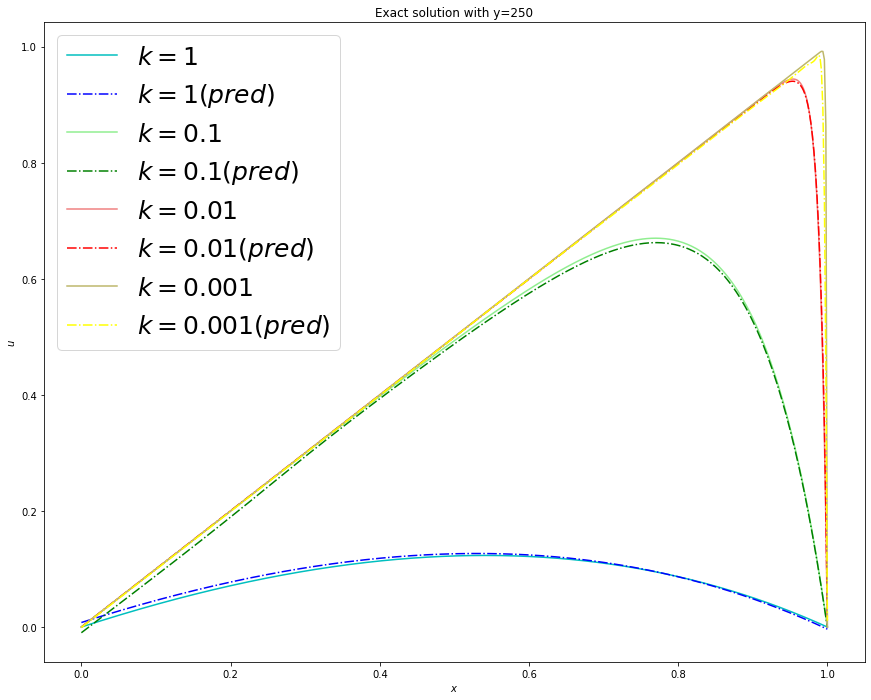}
    \caption{Scenario 2, advection-diffusion setting: predicted solution vs exact solution with respect to parameter $k$.}
    \label{fig:adv_results_scenario2}
\end{figure}

For Scenario 2, we add 20 different, randomly sampled $k \in (0.001, 1.0)$ to the input data.
In Figure~\ref{fig:adv_results_scenario2}, we show the exact solution and the solution field generated by a PINN with decreasing values of $k$ for this scenario.
The results are again impressive; the PINN algorithm in Scenario 2 is able to reconstruct the solution field with high precision from a small number of points used for the training, even for the case where we have a very small $k$. 

\begin{figure}[h!]
    \centering
    \includegraphics[width=.5\textwidth]{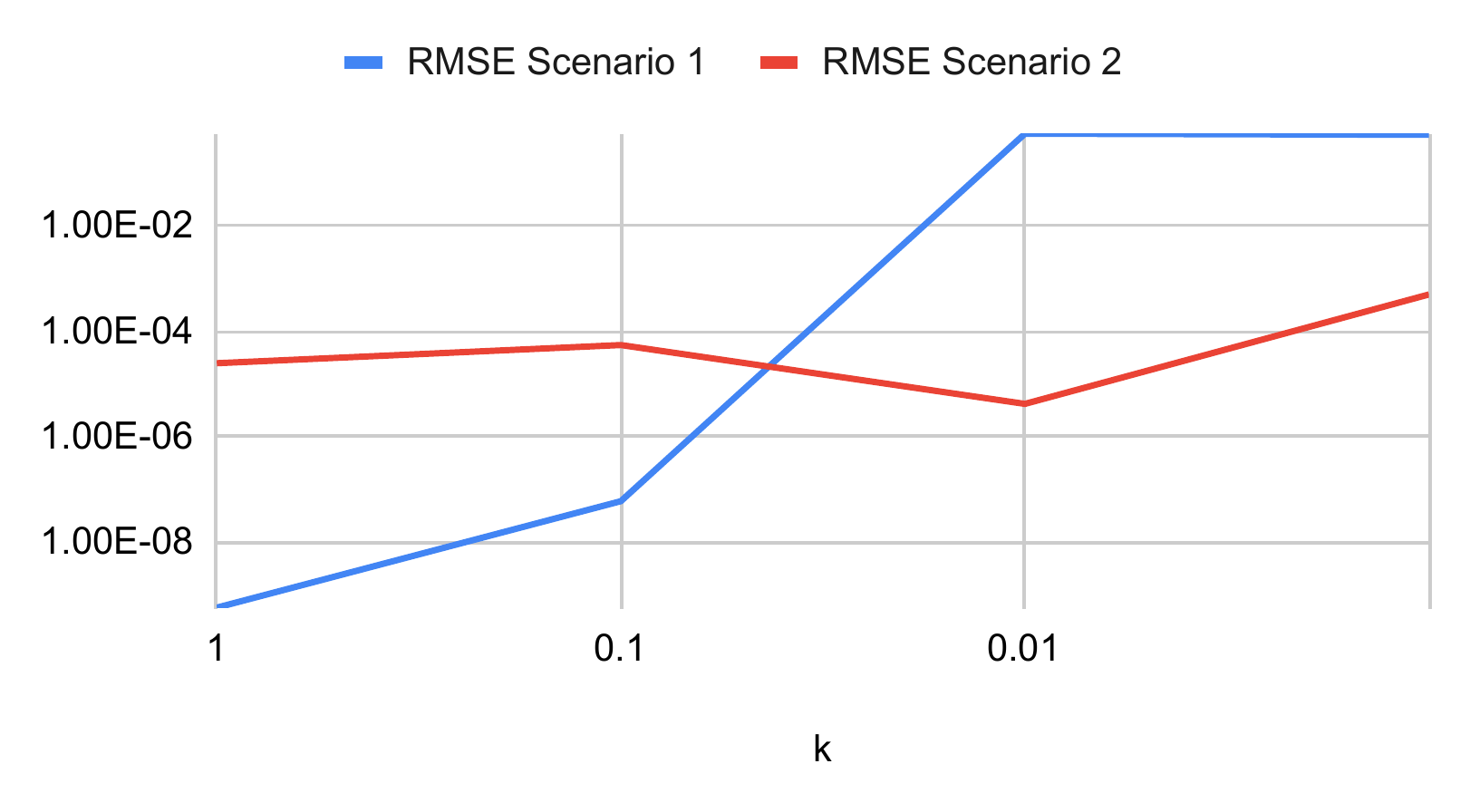}
    \caption{Quality of prediction for the advection dominant case (log-log scale).}
    \label{fig:advection_case}
\end{figure}

Lastly, Figure~\ref{fig:advection_case} presents the RMSE originating from Scenario 1 and Scenario 2. 
Again, for Scenario 2 the increase in the error is much slower than for Scenario 1.
Nevertheless, as in the reaction-diffusion case, the much higher errors for Scenario 2 with larger values of $k$ are still largely unexplained and motivates a series of investigations as part of our future work.

\section{Summary and Outlook}
\label{sec:conclusion}

The results presented herein clearly show the potential of PINNs for predicting the solution of PDEs with complex geometries and highly variable coefficients. 
Bringing physical coefficients into the training stage is key to avoiding the discrete solution's spurious oscillatory behavior in singularly perturbed regimes. 
In addition, it allows obtaining an accurate parameterization of the discrete solution concerning the physical coefficient in the interpolation scenario.
However, when we extrapolate, we can observe that this methodology will hardly overcome the numerical methods to solve direct linear problems and, while automatic hyperparameter selection is possible, it can be expensive.

One way to solve this problem is to combine the strengths of numerical methods and data science by creating hybrid combinations of theory-based and data science models. 
We will focus on the family of Multiscale Hybrid-Mixed (MHM) methods~\cite{Barrenechea2020poly,harder2015multiscale} and their interaction with PINNs. 
The MHM methods are attractive because of their approximation properties and massive parallelization capability, which allows the physical properties of the model to be treated locally and efficiently, thanks to the concept of local multiscale functions. 
So, we envision PINNs being used within MHM as surrogate models to predict the shape of the multiscale basis functions in parallel, among other possibilities. 
This combination will be the subject for future work.

Other topics for future work include: 
\begin{inparaenum}[(i)]
\item to explore alternatives to impose boundary conditions strongly (e.g.\ as in~\cite{EYu2018});
\item to apply the technique to other parametric studies, such as solutions with oscillatory behavior arising in oscillatory coefficient models or wave equation propagation problems;
\item to consider the use of PINNs for other expensive linear problems, such as in multi-query scenarios.
\end{inparaenum}

\section*{Acknowledgments}

The authors were supported by CAPES/Brazil under Project EOLIS (No.\ 88881.520197/2020-01).
Fr\'ed\'eric Valentin was supported by Inria/France under the Inria International Chair, by CNPq/Brazil under Project No.\ 309173/2020-5, and by FAPERJ/Brazil under Project No.\ E-26/201.182/2021.

\bibliography{pinn_boundarylayer}

\end{document}